\documentclass[11pt]{amsart}
\usepackage{amssymb, color}
\usepackage{amsmath}
\usepackage{amsthm}

\makeatletter
\def\LaTeX{\leavevmode L\raise.42ex
    \hbox{\kern-.3em\size{\sf@size}{0pt}\selectfont A}\kern-.15em\TeX}
\makeatother

\newcommand{\BibTeX}{{\rm B\kern-.05em{\sc
          i\kern-.025emb}\kern-.08em\TeX}}

\makeatletter
\def\@currentlabel{2.1}\label{e:dispaa}
\def\@currentlabel{2.21}\label{e:dispau}
\def\@currentlabel{2.22}\label{e:dispav}
\def\@currentlabel{2.23}\label{e:dispaw}
\def\@currentlabel{2.24}\label{e:dispax}
\def\theequation{\thesection.\@arabic\c@equation}
\makeatother

\usepackage{cases}
\newtheorem{theorem}{Theorem}[section]

\theoremstyle{definition}

\newtheorem{example}[theorem]{Example}
\theoremstyle{remark}

\numberwithin{equation}{section}
\newtheorem{conjecture}{Conjecture}[section]

\begin{document}
\title{Deformation of  $\mathfrak{osp}(2|2)$-Modules of Symbols}

\author{Mabrouk ben Ammar}
\address{Département de Mathématiques, Facult\'{e} des Sciences de Sfax, BP 802, 3038 Sfax, Tunisie}
\email{mabrouk.benammar@fss.rnu.tn}

\author{ Wafa Mtaouaa $^*$}
\address{Universit\'e de Sfax, Facult\'e des Sciences, D\'epartement de Math\'ematiques, B.P 1171 Sfax 3000, Tunisie.}
\email{mtawaa.wafa@yahoo.fr}



%

 \maketitle
%

\begin{abstract}
We classify deformations of $\mathfrak{osp}(2|2)-$module
structure on the spaces of symbols $\mathfrak{S}_d^2$ of
differential operators acting on the space of weighted densities
$\mathfrak{F}_{\lambda}^{2}$.\\
\textbf{Keywords:}  Cohomology, Orthosymplectic superalgebra,
deformation, weighted densities.  \\
\textbf{MSC(2010):} 17B56, 53D55, 58H15.
\end{abstract}

\section{\bf Introduction}

The space of weighted densities of weight $\mu$ on
$\mathbb{R}$ (or $\mu$-densities for short), denoted by:
\begin{equation*}{\mathcal F}_\mu=\left\{ fdx^{\mu}, ~f\in
C^\infty(\mathbb{R})\right\},\quad \mu\in\mathbb{R},
\end{equation*}is the space of
sections of the line bundle $(T^*\mathbb{R})^{\otimes^\mu}.$ The
Lie algebra ${\rm Vect}(\mathbb{R})$ of vector fields $X_h=h{d\over
dx}$, where $h\in C^\infty(\mathbb{R})$, acts by the {\it Lie
derivative}. Alternatively, this action can be written as follows:
\begin{equation}\label{Lie1}{X_h}\cdot(fdx^{\mu}):=L_{X_h}^\mu(fdx^{\mu})=(hf'+\mu h'f)dx^\mu,
\end{equation}
where $f'$ and $h'$ are $\frac{df}{dx}$ and $\frac{dh}{dx}$.

For $(\lambda,\mu)\in\mathbb{R}^2$ we consider  the space $\mathrm{D}_{\lambda,\mu}
:=\mathrm{Hom}_{\rm{diff}}(\mathcal{F}_{\lambda},\mathcal{F}_{\mu})$ of linear differential operators $A$ from $\mathcal{F}_\lambda$ to $\mathcal{F}_\mu$. The Lie algebra ${\rm Vect}(\mathbb{R})$ acts
on the space $\mathrm{D}_{\lambda,\mu}$ by:
\begin{equation}\label{Lieder2}X_h\cdot A=L_{X_h}^\mu\circ
A-A\circ L_{X_h}^{\lambda}. \end{equation}
Each module $\mathrm{D}_{\lambda,\mu}$ has a natural filtration by
the order of differential operators; the graded module $\mathcal
S_{\lambda,\mu}:=\mathrm{gr}\mathrm{D}_{\lambda,\mu}$ is called the
\textsl{space of symbols}. The quotient-module
$\mathrm{D}^k_{\lambda,\mu}/\mathrm{D}^{k-1}_{\lambda,\mu}$ is
isomorphic to $\mathcal
F_{\lambda-\mu-k}$, the isomorphism is provided by the principal
symbol $\sigma_{pr}$ defined by
$$
A=\sum_{i=0}^ka_i(x)\partial^i_x\mapsto\sigma_{pr}(A)=a_k(x)(dx)^{\mu-\lambda-k}$$
As a $\mathrm{Vect}(\mathbb{R})$-module, the space $\mathcal S_{\lambda,\mu}$
depends only on the difference $d=\mu-\lambda,$ so that
$\mathcal S_{\lambda,\mu}$ can be written as $\mathcal S_d$,
and we have
$$
\mathcal{S}_d=\bigoplus_{k=0}^{\infty}\mathcal{F}_{d-k}.$$
 Deformation problems appear in various areas of mathematics, in particular in algebra, algebraic and analytic geometry, and mathematical physics. Many powerful technics were developed to determine all related deformation obstructions.  The deformation theory of Lie algebras is widely studied. Some general questions of the theory were first considered by Richardson-Neijenhuis \cite{nr}. Their approach gave a strong relation between a given structure of Lie algebras and adapted cohomological tools. In fact, according to Richardson-Neijenhuis, deformation theory of modules is closely related to the computation of cohomology. In order to make this statement more precise, given a Lie algebra $\mathfrak{g}$ and a $\mathfrak{g}$-module $V$, the
infinitesimal deformations of the $\mathfrak{g}$-module structure on $V$, i.e., deformations that are linear in
the parameter of deformation are classified by the first cohomology space $\mathrm{H}^1(\mathfrak{g},\,\mathrm{End}(V))$. Of course, not for every infinitesimal deformation there exists a formal deformation
containing the latter as an infinitesimal part. The obstructions or conditions for which an infinitesimal deformation
guarantees existence of a formal deformation, are characterized
in terms of cup products (also called the Nijenhuis-Richardson products, see \cite{nr}) of non-trivial first cohomology
classes. These obstructions belong to the second cohomology
space $\mathrm{H}^2(\mathfrak{g},\,\mathrm{End}(V))$. This main result have been used by many authors (see \cite{aalo}, \cite{abbo}, \cite{bb2}, \cite{bbbbk}, \cite{bbdo}, \cite{bb} and references therein).

Consider the  $\mathrm{Vect}(\mathbb{R})$-module $\mathrm{D}_d:=\mathrm{Hom}_{\rm{diff}}(\mathcal{S}_d,\mathcal{S}_d)= \bigoplus_{i,j \geq0}\mathrm{D}_{d-i,d-j}$.  The space $\mathrm{D}_{\lambda,\mu}$ cannot be
isomorphic as a  $\mathrm{Vect}(\mathbb{R})$-module to the corresponding
space of symbols, but it is a deformation of this space in the sense of
Richardson-Neijenhuis \cite{nr}.

By restricting ourselves to the Lie algebra $\frak{ sl}(2)$ which
is isomorphic to the Lie subalgebra of ${\rm Vect}(\mathbb{R})$
spanned by
\begin{equation*}\left\{X_1,\,X_x,\,X_{x^2}\right\},\end{equation*}
we get families of infinite dimensional $\frak{ sl}(2)$-modules
still denoted by $\mathcal{F}_\lambda$, $\mathrm{D}_{\lambda,\mu}$ and $\mathcal{S}_d$.

Now, let us consider the superspace $\mathbb{R}^{1|n}$ endowed with its standard contact
structure defined by the 1-form $\alpha_n$, and the Lie superalgebra
$\mathcal{K}(n)$ of contact vector fields on
$\mathbb{R}^{1|n}$. We introduce the $\mathcal{K}(n)$-modules
$\mathfrak{F}_\lambda^n$ of $\lambda$-densities on
$\mathbb{R}^{1|n}$ and the $\mathcal{K}(n)$-modules of linear differential
operators, $\mathfrak{D}^n_{\lambda,\mu}
:=\mathrm{Hom}_{\rm{diff}}(\mathfrak{F}_{\lambda}^n,\mathfrak{F}_{\mu}^n)$, which are super
analogues of the spaces $\mathcal{F}_\lambda$ and
$\mathrm{D}_{\lambda,\mu}$, respectively. The module
$\frak{D}_{\lambda,\mu}^n$ is filtered:
\begin{equation*}
\mathfrak{D}^{n,0}_{\lambda,\mu}\subset\mathfrak{D}^{n,\frac{1}{2}}_{\lambda,\mu}\subset
\mathfrak{D}^{n,1}_{\lambda,\mu}\subset
\mathfrak{D}^{n,\frac{3}{2}}_{\lambda,\mu}
\subset\cdots\subset\mathfrak{D}^{n,\ell-\frac{1}{2}}_{\lambda,\mu}\subset\mathfrak{D}^{n,\ell}_{\lambda,\mu}
\cdots.
\end{equation*}
The corresponding graded module $\mathfrak
S_{\lambda,\mu}:=\mathrm{gr}\mathfrak{D}_{\lambda,\mu}^n$ is isomorphic to
 $$\mathfrak{S}^{n}_d = \bigoplus_{k=0}^\infty
\mathfrak{F}^{n}_{d-\frac{k}{2}},\quad d=\mu-\lambda.$$
We also consider the  $\mathcal{K}(n)$-module $\mathfrak{D}_d^n:=\mathrm{Hom}_{\rm{diff}}(\mathfrak{S}_d^n,\mathfrak{S}_d^n)= \bigoplus_{i,j \geq0}\mathfrak{D}_{d-\frac{i}{2},d-\frac{j}{2}}$.

The  Lie superalgebra $\mathfrak{osp}(n|2)$ is the super analogue of $\mathfrak{sl}(2)$ and it can be realized as a subalgebra of $\mathcal{K}(n)$. The spaces $\mathfrak{F}_\lambda^n$, $\mathfrak{D}^{n}_{\mu,\lambda}$ and $\mathfrak{S}^{n}_{\delta}$ are also $\mathfrak{osp}(n|2)$-modules.

We are interested to study the formal deformations of the
$\mathfrak{osp}(n|2)$-modules $\mathfrak{S}_d^{n}$.
According to Nijenhuis-Richardson \cite{nr}, the space
$\mathrm{H}^1\left(\mathfrak{osp}(n|2),\mathfrak{D}_d^{n}\right)$ classifies
the infinitesimal deformations of the $\mathfrak{osp}(n|2)$-module  $\mathfrak{S}_d^{n}$ and
the obstructions to integrability of a given infinitesimal
deformation of  $\mathfrak{S}_d^{n}$ are elements of
$\mathrm{H}^2\left(\mathfrak{osp}(2|n),\mathfrak{D}_d^{n}\right)$. For $n=0$ (classical case), the cohomology spaces $\mathrm{H}^1_\mathrm{diff}\left(\mathfrak{sl}(2),
\mathrm{D}_{\lambda,\mu}\right)$ and
$\mathrm{H}^2_\mathrm{diff}\left(\mathfrak{sl}(2),
\mathrm{D}_{\lambda,\mu}\right)$ were computed by Lecomte \cite{lec}. For $n=1$, Basdouri and Ben Ammar computed the cohomology spaces $\mathrm{H}^1_\mathrm{diff}\left(\mathfrak{osp}(1|2),
\mathfrak{D}_{\lambda,\mu}\right)$ \cite{bb1} and they studied the formal deformations of the $\mathfrak{sl}(2)$-modules $\mathcal{S}_d$ and the
$\mathfrak{osp}(1|2)$-modules $\mathfrak{S}_d^{1}$ \cite{bb2}. They exhibited
the necessary and sufficient integrability conditions of a given
infinitesimal deformation to a formal one and they proved that any formal deformation is
equivalent to its infinitesimal part. This work was generalized for $n\geq 3$ by Abdaoui, Khalfoun and Laraeidh \cite{akl} since in this case certain cohomological properties of the Lie superalgebras $\mathfrak{osp}(n|2)$ are similar. So, there seems to be no difference in results obtained in the study of  non-trivial deformations of the natural action of this orthosymplectic Lie superalgebra
on the direct sum of the superspaces of weighted densities. However, the case $n=2$ is exceptional because of an unexpected isomorphism $\mathcal{K}(n))\simeq\mathrm{Vect}(\mathbb{R}^{1|1})$ (see \cite{gro}) which motivate Ben Fraj and Boujelben in \cite{bb3} to compute the cohomology space $\mathrm{H}^{1}_{\mathrm{diff}}(\mathfrak{osp}(2|2),\mathfrak{D}^{2}_{\lambda\,\mu})$.

In this paper we are interested in the case $n=2$, we study the formal deformations of the
$\mathfrak{osp}(2|2)$-modules $\mathfrak{S}_d^{2}$ and we give the necessary and sufficient integrability conditions of a given
infinitesimal deformation to a formal one.

\section{Definitions and Notations}

Let $\mathbb{R}^{1|2}$ be the superspace with coordinates $(x,\,\theta_1,\,\theta_2)$ where
$x$ is the even indeterminate, $\theta_1$ and $\theta_2$ are odd
 indeterminates, i.e., $\theta_i\theta_j=-\theta_j\theta_i$. $\mathbb{R}^{1|2}$ is equipped with the standard contact structure given by the following $1$-form:
\begin{equation*}
\alpha_2=dx+\theta_1 d\theta_1+\theta_2 d\theta_2.
\end{equation*}
Consider $C^\infty({\mathbb{R}}^{1|2})$ through the space of
functions of $\mathbb{R}^{1|2}$. A $C^\infty({\mathbb{R}}^{1|2})$ function has the form:
\begin{equation*}
f(x,\theta_1,\theta_2)=f_0(x) + f_1(x)\theta_1 + f_2(x)\theta_2 +
f_{12}(x)\theta_1\theta_2,
\end{equation*}
where $f_0, f_1, f_2, f_{12}\in C^\infty({\mathbb{R}})$. We denote by $|f|$
the parity of an homogeneous function $f$, that is, $|f_0(x)|=|f_{12}(x)\theta_1\theta_2|=0$ and $|f_1(x)\theta_1|=|f_2(x)\theta_2|=1$. Hereafter, the expression $(-1)^{|f|}$ will be simply written $(-1)^{f}$.

Let $\mathrm{Vect}(\mathbb{R}^{1|2})$ the superspace of vector fields on $\mathbb{R}^{1|2}$:
\begin{equation*}
\mathrm{Vect}(\mathbb{R}^{1|2})=\left\{h_0\partial_x+h_1\partial_{1}+h_2\partial_{2}\mid
h_i\in C^\infty(\mathbb{R}^{1|2})\right\},
\end{equation*}
where $\partial_x=\frac{\partial}{\partial x}$ and $\partial_{i}=\frac{\partial}{\partial\theta_i}$, and consider the Lie superalgebra  $\mathcal{K}(2)$ of contact vector fields on
$\mathbb{R}^{1|2}$. That is, $\mathcal{K}(2)$ is a subalgebra of $\mathrm{Vect}(\mathbb{R}^{1|2})$ preserving the
distribution singled out by the 1-form $\alpha_2$:
$$\mathcal{K}(2)=\big\{X\in\mathrm{Vect}(\mathbb{R}^{1|2})~|~\hbox{there
exists}~f\in C^\infty({\mathbb{R}}^{1|2})~ \hbox{such
that}~\mathfrak{L}_X(\alpha_2)=f\alpha_2\big\},
$$
where $\mathfrak{L}_X$ is the Lie derivative along the vector field
$X$.\\
 Consider the vector fields $\overline{\eta}_i=\partial_{i}-\theta_i\partial_x$, any contact vector field on $\mathbb{R}^{1|2}$ can be
expressed as
\begin{equation*}
X_f=f\partial_x-\frac{1}{2}(-1)^{f}\sum_{i=1}^2\overline{\eta}_i(f)\overline{\eta}_i,\;\text{
where }\, f\in C^\infty(\mathbb{R}^{1|2}).
\end{equation*}
 The contact
bracket is defined by $[X_f,\,X_g]=X_{\{f,\,g\}}$ where $\{,\,\}$ is the Poisson bracket defined by
\begin{equation}
\{f,\,g\}=fg'-f'g-\frac{1}{2}(-1)^{f}\sum_{i=1}^2\overline{\eta}_i(f)\cdot\overline{\eta}_i(g).
\end{equation}
Then the map $f\mapsto X_f$ is an isomorphism of Lie superalgebra from $( C^\infty(\mathbb{R}^{1|2}),\,\{,\,\})$ to $(\mathcal{K}(2),\,[,\,])$.
Thus, via this isomorphism, the Lie superalgebra $\mathcal{K}(2)$ can be identified to the Lie superalgebra  $C^\infty(\mathbb{R}^{1|2})$ endowed with the Poisson bracket.

We define the Lie superalgebra
\begin{equation*}\mathfrak{osp}(2|2)=\langle H,X,Y,A_1,A_2,B_1,B_2,C\rangle.\end{equation*}
 The elements $H$, $X$, $Y$ and $C$ are even and the elements $A_i$, $B_i$ are odd, the bracket is graded antisymmetric, we denote this
 property by $$[U,V]=-(-1)^{UV}[V,U].$$ The non zero brackets are:
$$\begin{array}{llll}
&[A_i,A_i]=2X,~~&[X,Y]=2H,~~&[H,X]=X,\\
&[A_i,Y]=-B_i,~~&[X,B_i]=A_i,~~&[H,A_i]=\frac{1}{2}A_i,\\
&[A_i,B_i]=2H,~~&[B_i,B_i]=-2Y,~~&[H,B_i]=-\frac{1}{2}B_i,\\
&[A_1,C]=\frac{1}{2}A_2,~~&[B_1,C]=\frac{1}{2}B_2,~~&[H,Y]=-Y,\\
&[A_2,C]=-\frac{1}{2}A_1,~~&[B_2,C]=-\frac{1}{2}B_1.
\end{array}
$$
It is well known that $\mathfrak{osp}(2|2)$ can be realized as a subalgebra of $\mathcal{K}(2)$:
\begin{equation*}\mathfrak{osp}(2|2)=\text{Span}\left(1,\,{x},\,{x^2},\,{x\theta_1},\,
{x\theta_2},\, {\theta_1},\, {\theta_2},\,
\theta_1\theta_2\right),\end{equation*}
 Here,
\begin{equation*}\left(-x,\,1,\,-{x^2},\,2{\theta_i},\,2{x\theta_i},\,\theta_1\theta_2\right)=(H,\,X,\,Y,\,A_i,\,B_i,\,C)
\end{equation*}
We easily see that
$\mathfrak{osp}(1|2)$ is isomorphic to a subalgebra of $\mathfrak{osp}(2|2)$:
$$
\mathfrak{osp}(1|2)\simeq\mathfrak{osp}(1|2)^i=
\text{Span}\left(1,\,{x},\,{x^2},\,{x\theta_i},\, {\theta_i}\right),\quad i=1,\,2.
$$
We define the space of $\lambda$-densities as
\begin{equation}
\mathfrak{F}^2_\lambda=\left\{f(x,\theta_1,\theta_2)\alpha_2^\lambda\mid
f(x,\theta_1,\theta_2) \in C^\infty(\mathbb{R}^{1|2})\right\}.
\end{equation}
As a vector space, $\mathfrak{F}^2_\lambda$ is isomorphic to
$C^\infty(\mathbb{R}^{1|2})$, but the Lie derivative of the density
$g\alpha_2^\lambda$ along the vector field $f:=X_f$ in $\mathcal{K}(2)$
is now:
\begin{equation}
\label{superaction}
\mathfrak{L}^{\lambda}_{f}(g\alpha_2^\lambda)=(\mathfrak{L}_{f}(g)+
\lambda f'g)\alpha_2^\lambda.
\end{equation}

Here, we restrict ourselves to the subalgebra $\mathfrak{osp}(2|2)$,
thus we obtain a one-parameter family of
$\mathfrak{osp}(2|2)$-modules on $C^\infty(\mathbb{R}^{1|2})$ still
denoted by $\mathfrak{F}^2_\lambda$. As an
$\mathfrak{osp}(1|2)$-module, we have
\begin{equation}\label{isom}
\mathfrak{F}^2_\lambda\simeq\mathfrak{F}^1_\lambda
\oplus\Pi(\mathfrak{F}^1_{\lambda+{1\over2}})
\end{equation}
where $\Pi$ is the change of parity operator.
\section{Cohomology}
 Let $\frak{g}$ be a Lie superalgebra acting on a
superspace $V$.  The space of $n$-cochains of $\frak{g}$ with values in $V$ is the
$\frak{g}$-module
\begin{equation*}
C^n(\frak{g}, V ) := \mathrm{Hom}(\Lambda^n(\frak{g}),V).
\end{equation*}
The {\it coboundary operator} $ \delta_n: C^n(\frak{g}, V
)\longrightarrow C^{n+1}(\frak{g}, V )$ is a
$\frak{g}$-map satisfying $\delta_n\circ\delta_{n-1}=0$. The
kernel of $\delta_n$, denoted $Z^n(\mathfrak{g},V)$, is
the space of $n$-{\it cocycles}, among them,
the elements in the range of $\delta_{n-1}$ are called $n$-{\it coboundaries}. We denote
$B^n(\mathfrak{g},V)$ the space of $n$-coboundaries.
By definition, the $n^{th}$   cohomology space is
the quotient space
\begin{equation*}
\mathrm{H}^n
(\mathfrak{g},V)=Z^n(\mathfrak{g},V)/B^n(\mathfrak{g},V).
\end{equation*}
We will only need the formula of $\delta_n$ (which will be simply
denoted $\delta$) in degrees 0 and 1: for $v \in
C^0(\frak{g}, V) =V$,~ $\delta v(g) : =
(-1)^{gv}g\cdot v$
and  for  $ \omega\in C^1(\frak{g}, V )$,
\begin{equation*}\delta \omega(g,\,h):=
(-1)^{g\omega}g\cdot
\omega(h)-(-1)^{h(g+\omega)}h\cdot
\omega(g)-\omega([g,~h])\quad\text{for any}\quad g,h\in
\frak{g}.
\end{equation*}
For the general expression of $\delta_n$ see eg \cite{ddar}.

\section{Deformation theory and cohomology}
 Let $\rho_{0}:\mathfrak{g}\longrightarrow \mathrm{End}(V)$ be an action of a
Lie superalgebra $\mathfrak{g}$ on a vector superspace $V$. When studying
deformations of the $\mathfrak{g}$-action $\rho_{0}$, one usually starts with
infinitesimal deformations:
\begin{equation}
\rho=\rho_{0}+t\omega
\end{equation}
where $\omega:\mathfrak{g}\rightarrow \mathrm{End}(V)$ is a linear map and
$t$ is a formal parameter with $|t|=|\omega|$. From the homomorphism
condition
\begin{equation}
[\rho(x)\,,\rho(y)]=\rho([x,y])
\end{equation}
where $x,\,y \in \mathfrak{g}$, we deduce that
$\omega$ is an $1$-cocycle. That is, the linear map $\omega$ satisfies
\begin{equation}
(-1)^{x\omega}[\rho_{0}(x),\omega(y)]-(-1)^{y(x+\omega)}[\rho_{0}(y),\omega(x)]- \omega([x,y])=0.
\end{equation}
Moreover, two infinitesimal deformations
$\rho=\rho_{0}+t\omega_{1}$, and $\rho=\rho_{0}+t\omega_{2}$ are
equivalents if and only if $c_{1}-c_{2}$ is coboundary:
\begin{equation}
(\omega_{1}-\omega_{2})(x)=(-1)^{xA}[\rho_{0},\,A](x):=\delta A(x)
\end{equation}
where $A\in \mathrm{End}(V)$. So,
the space $\mathrm{H}^{1}(\mathfrak{g},\,V)$ determines and classifies
infinitesimal deformations up to equivalence.

Now, if $\dim
(\mathrm{H}^{1}(\mathfrak{g},\,V))=m$, then we choose $1$-cocycles
$\omega_{1},\dots, \omega_{m}$ representing a basis of
$\mathrm{H}^{1}(\mathfrak{g},\,V)$ and we consider the infinitesimal deformation
\begin{equation}
\rho=\rho_{0}+\sum\limits_{i=1}^{m}t_{i}\omega_{i}
\end{equation}
where $t_{1},\dots,t_{m}$ are independent parameters with
$|t_{i}|=|\omega_{i}|$. We try to extend this infinitesimal
deformation to a formal one
\begin{equation}
\rho=\rho_{0}+\sum\limits_{i=1}^{m}t_{i}\omega_{i}+\sum\limits_{i,j}t_{i}t_{j}\rho_{ij}^{2}+\cdots
\end{equation}
where $\rho_{ij}^{2},\,\rho_{ijk}^{3}\cdots$ are linear maps from
$\mathfrak{g}$ to $\mathrm{End}(V)$ with
$|\rho_{ij}^{2}|=|t_{i}t_{j}|,\,|\rho_{ijk}^{3}|=|t_{i}t_{j}t_{k}|$
such that
\begin{equation}
[\rho(x),\,\rho(y)]=\rho([x,\,y]),\quad x,\,y \in \mathfrak{g}
\end{equation}
All the obstructions become from this condition  and it is well known that they lie in $\mathrm{H}^{2}(\mathfrak{g},\,V)$. Thus, we will
impose extra algebraic relations on the parameters $t_{1},\dots,
t_{m}$.
Let $\mathcal{R}$ be an an ideal in $\mathbb{C}[[t_{1},\dots,
t_{m}]]$ generated by some set of these relations, the quotient
\begin{equation}
\mathcal{A}=\mathbb{C}[[t_{1},\dots, t_{m}]]/\mathcal{R}
\end{equation}
is a supercommutative associative superalgebra with unity.
\section{Cohomology and deformation of $\mathfrak{S}_d^{2}$}
We study the formal deformations of the $\mathfrak{osp}(2|2)$-module
structure on the space of symbols:
$$
\mathfrak{S}_d^{2}
=\bigoplus\limits_{k\geq0}\mathfrak{F}^{2}_{d-\frac{k}{2}}
$$
The infinitesimal deformations are described by the cohomology
space:
  $$\mathrm{H}^{1}_{\mathrm{diff}}\left(\mathfrak{osp}(2|2),\mathfrak{S}_d^{2}\right)=\bigoplus\limits _{i,j\geq0}\mathrm{H}^{1}_{\mathrm{diff}}\left(\mathfrak{osp}(2|2), \mathfrak{D}_{d-\frac{j}{2},d-\frac{i}{2}}^{2}\right)$$

  Ben Fraj and Boujelben \cite{bb3} computed the spaces
 $\mathrm{H}^{1}_{\mathrm{diff}}(\mathfrak{osp}(2|2),\mathfrak{D}_{\lambda,\mu}^{2})$, they showed the following result:
\begin{theorem}  \begin{equation} \dim(\mathrm{H^1_{diff}}(\mathfrak{osp}(2|2),\frak{D}^2_{\lambda,\mu}))=\left\{ \begin{array}{llllll} 2&\text{if}\quad
\lambda=\mu,\\[2pt] 3&\text{if}\quad (\lambda,\mu)=(-\frac{k}{2},\frac{k}{2}) \hbox{
with } k\in\mathbb{N}\backslash\{ 0\},\\[2pt] 0&\text{otherwise}. \end{array} \right. \end{equation}
\end{theorem}
Moreover,  basis for these  cohomology spaces are given   in \cite{bb3}.
Thus,
\begin{itemize}
  \item [i)]
 If $2d \notin \mathbb{N}$, then
 $$\mathrm{H}^{1}_{\mathrm{diff}}\left(\mathfrak{osp}(2|2),\mathfrak{S}_d^{2}\right)= \bigoplus\limits_{k\geq 0}\mathrm{H}^{1}_{\mathrm{diff}}\left(\mathfrak{osp}(2|2), \mathfrak{D}_{d-\frac{k}{2},d-\frac{k}{2}}^{2}\right).$$
 The space
 $\mathrm{H}^{1}_{\mathrm{diff}}\left(\mathfrak{osp}(2|2), \mathfrak{D}_{d-\frac{k}{2},d-\frac{k}{2}}\right)$
  is spanned by:
$$\omega_{k}({f})=f'\quad\text{ and }\quad
\widetilde{\omega}_{k}({f})=(2d-k)\overline{\eta}_{1}\partial_{2}f- (-1)^{f}(\partial_{2}f\overline{\eta}_{1}+ \theta_{2}\overline{\eta}_{2}\overline{\eta}_{1}f\overline{\eta}_{2}).$$
 \item [ii)]
 If $ 2d=m \in \mathbb{N}$, then
 $$\mathrm{H}^{1}_{\mathrm{diff}}\left(\mathfrak{osp}(2|2),\mathfrak{S}_d^{2}\right)= \bigoplus\limits_{k=1}^{m}\mathrm{H}^{1}_{\mathrm{diff}}\left(\mathfrak{osp}(2|2), \mathfrak{D}_{-\frac{k}{2},\frac{k}{2}}^{2}\right)
 \oplus \bigoplus\limits_{k=-\infty}^m\mathrm{H}^{1}_{\mathrm{diff}}\left(\mathfrak{osp}(2|2), \mathfrak{D}_{\frac{k}{2},\frac{k}{2}}^{2}\right).$$
  The space
 $\mathrm{H}^{1}_{\mathrm{diff}}\left(\mathfrak{osp}(2|2),\mathfrak{D}_{\frac{k}{2},\frac{k}{2}}\right)$
 is  spanned by:
 \begin{equation}\label{maincocyc}
\begin{array}{lllllllllll}\gamma_{k}({f})&=&f'\\[5pt]
\widetilde{\gamma}_{k}({f})&=&\left\{\begin{array}{ll}
\overline{\eta}_1\overline{\eta}_2f\hfill\text{
 if }k=0\\[2pt]k\,
\overline{\eta}_1\partial_{2}f-(-1)^{f}
\left(\partial_{2}f\overline{\eta}_1+\theta_2\overline{\eta}_2
\overline{\eta}_{1}f\overline{\eta}_2 \right)\hfill\text{ if
}k\neq 0.\end{array}\right. \end{array}\end{equation}
 The space
 $\mathrm{H}^{1}_{\mathrm{diff}}\left(\mathfrak{osp}(2|2),\mathfrak{D}_{-\frac{k}{2},\frac{k}{2}}^{2}\right)$
  is spanned by:
 \begin{eqnarray*}
 \Gamma_{k}({f})&=&f'\overline{\eta}_{1}\overline{\eta}_{2}^{2k-1}\\
 \widetilde{\Gamma}_{k}({f})&=&k\overline{\eta}_{1}(\partial_{2}f)\overline{\eta}_{1} \overline{\eta}_{2}^{2k-1}-(-1)^{f}\left(\partial_{2}f\overline{\eta}_{2}^{2k+1} -\overline{\eta}_{1}(\theta_{2}\partial_{2}f)\overline{\eta}_{1}^{2k+1}\right)\\
 \overline{\Gamma}_{k}({f})&=&(k-1)f''\overline{\eta}_{1}\overline{\eta}_{2}^{2k-3}+
 (-1)^{f}\left(\overline{\eta}_{2}f'\overline{\eta}_{1}^{2k-1}- \overline{\eta}_{1}f'\overline{\eta}_{2}^{2k-1}\right).
\end{eqnarray*}
\end{itemize}
 In our study,  any  infinitesimal deformation of $\mathfrak{osp}(2|2)$-module on the space
 $\mathfrak{S}_d^{2}$
 is of the form: \begin{equation}\label{inf}\widetilde{\mathfrak{L}}=\mathfrak{L}+ \mathfrak{L}^{1} \end{equation}
 where
 $$\mathfrak{L}^{1}= \left\{
\begin{array}{ll}
\sum_{k\geq 0} (a_{k}\omega_{k}+b_{k}\widetilde{\omega}_{k})&\text{if}\quad 2d\notin\mathbb{N}\\
\sum _{k\leq
m}(a_{k}\gamma_{k}+b_{k}\widetilde{\gamma}_{k})+\sum
_{k=1}^{m}(c_{k}\Gamma_{k}+d_{k}\widetilde{\Gamma}_{k}+e_{k}\overline{\Gamma}_k) &
\text{if}\quad
2d=m \in \mathbb{N}.
 \end{array}
 \right.$$
 The coefficients $a_{k},b_{k},c_{k},d_{k},e_{k}$ are  independent parameters.

Now, we extend the
infinitesimal deformation \eqref{inf} to a formal
one:
\begin{equation}
\label{BigDef2} \widetilde{\mathfrak{L}}= \mathfrak{L}+\mathfrak{L}^{1}+\sum_iP_i^{2}\mathfrak{L}_i^{2}+ \sum_i P_i^{3}\mathfrak{L}_i^{3}+\cdots,
\end{equation}
where the higher order terms
$\mathfrak{L}_i^{2}$, $\mathfrak{L}_i^{3},\ldots$ are linear maps from
$\mathfrak{osp}(2|2)$ to $\mathrm{End}(\mathfrak{S}_d^{2})$ such that the map
\begin{equation} \label{map} \widetilde{\mathfrak{L}}:\mathfrak{osp}(2|2)\to
\mathbb{C}[[a_{k},b_{k},c_{k},d_{k},e_{k}]]\otimes{\rm End(\mathfrak{S}_d^{2})},
\end{equation}
satisfies the homomorphism condition
 \begin{equation}\label{hom}\widetilde{\mathfrak{L}}_{[f,g]}= [\widetilde{\mathfrak{L}}_f,\widetilde{\mathfrak{L}}_g].\end{equation}
 $P_i^{j}$ are monomial in the parameters $a_{k},b_{k},c_{k},d_{k},e_{k}$ (or $a_{k},b_{k}$ if $2d\notin\mathbb{N}$) with degree $j$ and with the same parity of $\mathfrak{L}_i^{j}$.\\
Setting
\begin{equation*}\varphi = \widetilde{\mathfrak{L}}- \mathfrak{L},\quad
\mathfrak{L}^{2}=\sum_iP_i^{2}\mathfrak{L}_i^{2},\quad \mathfrak{L}^{3}=\sum_i P_i^{3}\mathfrak{L}_i^{3},\ldots,
\end{equation*}
we can rewrite the homomorphism condition \eqref{hom} in the following way:
\begin{equation}
\label{developping} [\varphi(f) , \mathfrak{L}_g ] + [\mathfrak{L}_f ,
\varphi(g) ] - \varphi([f , g]) +\sum_{i,j > 0}
\;[\mathfrak{L}^{i}_f , \mathfrak{L}^{j}_g] = 0,
\end{equation}
or equivalently
\begin{equation}
\label{maurrer cartan} \delta\varphi +{1\over2} \varphi \vee
\varphi= 0,
\end{equation}
where $\delta\varphi$ stands for differential
of the cochain $\varphi$ and $\vee$ is the standard {\it
cup-product} defined, for arbitrary
linear maps $a,~b :\mathfrak{g}
\longrightarrow\mathrm{End}(V)$ with $\mathfrak{g}$ a Lie superalgebra and $V$ a vector superspace,  by:
\begin{equation} \label{maurrer} (a\vee b) (x , y) =
(-1)^{xb}[a(x) , b(y)] +
(-1)^{a(x+b)}[b(x) , a(y)],
\end{equation}
so that, if $ a$ and $b$ are even maps then
\begin{equation*}
\label{maurrer cartan1} (a \vee b) (x , y) =
[a(x) , b(y)] +
[b(x) , a(y)].
\end{equation*}
From (\ref{maurrer cartan})  we obtain the following equation for any $\mathfrak{L}^{k}$:
\begin{equation}
\label{maurrer cartank} \delta\mathfrak{L}^{k} + {1\over2}\sum_{i+j=k}
\mathfrak{L}^{i} \vee  \mathfrak{L}^{j}= 0.
\end{equation}
The first non-trivial relation $$\delta{\mathfrak{L}^{2}} +{1\over2}
\mathfrak{L}^{1} \vee \mathfrak{L}^{1} = 0 $$ gives the first obstruction
to integration of an infinitesimal deformation. That is, $\mathfrak{L}^{1} \vee \mathfrak{L}^{1}$ must be a a coboundary.

It is easy to check that for any two $1$-cocycles $C_1$ and
$C_2 \in Z^1 (\frak g , \mathrm{End}(V))$, the
bilinear map $C_1 \vee C_2$ is a
$2$-cocycle.  Moreover, if one of the
cocycles $C_1$ or $C_2$ is a
coboundary, then $C_1 \vee C_2$ is a  $2$-coboundary. Therefore, we naturally deduce that
the operation (\ref{maurrer}) defines a bilinear map:
\begin{equation}
\label{cup-product} \mathrm{H}^1 (\frak g ,\mathrm{End}(
V))\otimes \mathrm{H}^1 (\frak g , \mathrm{End}(
V))\longrightarrow \mathrm{H}^2 (\frak g , \mathrm{End}(
V)).
\end{equation}
All the obstructions lie in $\mathrm{H}^2 (\frak g,\mathrm{End}(V))$ and they are in the image of $\mathrm{H}^1
(\frak g,\mathrm{End}(V))$ under the cup-product. Thus, we describe in the following section the cup-product $\mathrm{H}^1\vee\mathrm{H}^1$.

\section{The cup-product $\mathrm{H}^1\vee\mathrm{H}^1$}
We have to distinguish two cases:

\subsection{Case 1: $2d \notin \mathbb{N}$}

\begin{theorem} \label{th1}
If $2d \notin \mathbb{N}$ then the image $\mathrm{H}^1\vee\mathrm{H}^1$ of $\mathrm{H}^{1}_{\mathrm{diff}}\left(\mathfrak{osp}(2|2), \mathfrak{D}_{d-\frac{k}{2},d-\frac{k}{2}}\right)$ under the cup-product is a 2-dimensional subspace of $\mathrm{H}^{2}_{\mathrm{diff}}\left(\mathfrak{osp}(2|2), \mathfrak{D}_{d-\frac{k}{2},d-\frac{k}{2}}\right)$ spanned by $$\Omega_1=\omega_k\vee\widetilde{\omega}_k\quad\text{ and }\quad \Omega_2=\widetilde{\omega}_k\vee\widetilde{\omega}_k.$$
\end{theorem}
\begin{proofname}. In this case, the space $\mathrm{H}^1\vee\mathrm{H}^1$ is generated by the three cup-products: $\omega_k\vee\omega_k$, $\omega_k\vee\widetilde{\omega}_k$ and $\widetilde{\omega}_k\vee\widetilde{\omega}_k$. But it is easily check that $\omega_k\vee\omega_k=0$. So, we have to prove that $\Omega_1$ and $\Omega_2$ are nontrivial 2-cocycles which are linearly independent. That is, the equation
\begin{equation}\label{eq1}
a\Omega_1+b\Omega_2=\delta B,
\end{equation}
where $a,\,b\in\mathbb{R}$ and
$B\in \mathrm{C}^{1}_{\mathrm{diff}}\left(\mathfrak{osp}(2|2), \mathfrak{D}_{d-\frac{k}{2},
d-\frac{k}{2}}\right)$, has a solution if and only if $a=b=0$.

First of all, we have
 \begin{eqnarray*}
  \Omega_1(g,h)   &=&-(-1)^{g}\partial_2g\overline{\eta}_{1}h'
   -(-1)^{g}\theta_2\overline{\eta}_{2}\overline{\eta}_{1}g\overline{\eta}_{2}h'
   -(-1)^{gh}(g\leftrightarrow h), \\
   \Omega_2(g,h) &=&2\bigg[(-1)^{g}(2d-k)\partial_2g\partial_x\partial_2h+(-1)^{g}\partial_2g\textbf{(}\partial_2h-\theta_1\partial_1\partial_2h\textbf{)}\partial_x \\&&+ (-1)^{g+h} \partial_2g\overline{\eta}_{1}\partial_2h\partial_1 +
   (-1)^{h} \theta_2  \partial_2g\partial_x\partial_2h \partial_2\bigg]
  -(-1)^{gh}(g\leftrightarrow h). \end{eqnarray*}

Now, for $\alpha=(i,j,k)$, we denote by $\partial^\alpha=\partial_x^i\partial_1^{j}\partial_2^{k}$. Then, by considering the equation \eqref{eq1}, we can write
\begin{align}\label{cobo}
B (h)=\sum\limits_{\alpha,\beta}A_{\alpha,\beta}\partial^{\alpha}(h)\partial^{\beta}\quad \text{where}\quad A_{\alpha,\beta}= A^{0}_{\alpha,\beta}+\theta_1A^{1}_{\alpha,\beta}+\theta_2A^{2}_{\alpha,\beta}+\theta_1\theta_2A^{12}_{\alpha,\beta}.
  \end{align}
   One obtains
 \begin{eqnarray*}
  B (1)&=&\sum\limits_{\beta}A_{000,\beta}\partial^{\beta}\\
B (x)&=&\sum\limits_{\beta}(A_{000,\beta}x+A_{100,\beta})\partial^{\beta}\\
B (x^{2})&=&\sum\limits_{\beta}(A_{000,\beta}x^2+2xA_{100,\beta}+2A_{200,\beta})\partial^{\beta}\\
B(\theta_1)&=&\sum\limits_{\beta}(A_{000,\beta}\theta_1+A_{010,\beta})\partial^{\beta}\\
B(\theta_2)&=&\sum\limits_{\beta}(A_{000,\beta}\theta_2+A_{001,\beta})\partial^{\beta}\\
B(x\theta_1)&=&\sum\limits_{\beta}( A_{000,\beta}x \theta_1+ A_{100,\beta}\theta_1+ A_{010,\beta}x+A_{110,\beta}                   )\partial^{\beta}\\
B(x\theta_2)&=&\sum\limits_{\beta}( A_{000,\beta}x \theta_2+ A_{100,\beta}\theta_2+ A_{001,\beta}x+A_{101,\beta}                   )\partial^{\beta}\\
B(\theta_1\theta_2)&=&\sum\limits_{\beta}( A_{000,\beta}\theta_1\theta_2+ A_{010,\beta}\theta_2- A_{001,\beta}\theta_1+ A_{011,\beta}                 )\partial^{\beta}
 \end{eqnarray*}
Let us recall that
 \begin{eqnarray}
\delta B(g,\,h): &=& \mathfrak{L}_{g}^{\lambda,\mu}B(h) -(-1)^{hg} \mathfrak{L}_{g}^{\lambda,\mu}B(g)-B([g,h]) \nonumber\\
 &=&
 \partial_x B(h)-\frac{1}{2}(-1)^{g}\big( \overline{\eta}_1g \overline{\eta}_1 B(h)+\overline{\eta}_2g \overline{\eta}_2 B(h)\big)\nonumber \\
 && +\mu \partial_x g B(h)-B(h)\big(g \partial_x - \frac{1}{2}(-1)^{g} \big( \overline{\eta}_1g \overline{\eta}_1+\overline{\eta}_2g \overline{\eta}_2 \big)+\lambda \partial_x g \big) \nonumber\\
 &&-(-1)^{gh} \big( h\partial_x B(g)-\frac{1}{2}(-1)^{h}\big( \overline{\eta}_1h \overline{\eta}_1 B(g)+\overline{\eta}_2h \overline{\eta}_2 B(g)\big)\nonumber \\
 && +\mu \partial_xh B(g)-B(g)\big(h \partial_x - \frac{1}{2}(-1)^{h} \big( \overline{\eta}_1h \overline{\eta}_1+\overline{\eta}_2h \overline{\eta}_2 \big)+\lambda \partial_x h \big)  \big)\nonumber\\
 &&-B \big(g\partial_x h-\partial_xgh- \frac{1}{2}(-1)^{g}\big( \overline{\eta}_1g \overline{\eta}_1 h+\overline{\eta}_2g \overline{\eta}_2 h\big)\big)\nonumber
\end{eqnarray}

Now, considering the terms in $f$ in \eqref{eq1} for $(g,\,h)=(\theta_2,\theta_2)$ then for $(g,\,h)=(x\theta_2,\theta_2)$, we get
\begin{align}\label{a}
 -\lambda A^{0}_{001,001}+\frac{1}{4}A^{2}_{101,000}=-4\lambda b.
\end{align}
Similarly, the terms in ${\theta_1f}$ for $(g,\,h)=(\theta_1,\theta_1)$ then for $(g,\,h)=(x\theta_1,\theta_1)$ give
\begin{align}\label{b}
\frac{1}{4}A^{1}_{110,000}-\lambda  A^{0}_{010,010} =0.
\end{align}
The terms in ${\theta_2 f}$ for $(g,\,h)=(\theta_2,\theta_1\theta_2)$ then for $(g,\,h)=(x\theta_2,\theta_1\theta_2)$ give
 \begin{align}\label{e}
\lambda A^{0}_{001,001}- \lambda A^{1}_{011,001}-\frac{1}{4}A^{0}_{100,000}-\frac{1}{2}A^{2}_{101,000}+\frac{1}{4}A^{1}_{110,000}=4b\lambda.
\end{align}
Considering the terms in ${\partial_1 f}$ for $(g,\,h)=(\theta_1,\theta_1\theta_2)$ and $(g,\,h)=(x\theta_1,\theta_1)$ we obtain
 \begin{align}\label{f}
-\lambda A^{2}_{011,010}- \lambda A^{0}_{010,010}+\frac{1}{2}A^{1}_{110,000}-\frac{1}{4}A^{2}_{101,000}+\frac{1}{4}A^{0}_{100,000}=0.
\end{align}
Now, we consider the terms in ${\partial_2 f}$ respectively for $(g,h)=(\theta_2,\theta_1\theta_2)$ and for $(g,h)=(\theta_1,\theta_1\theta_2)$ then we obtain
\begin{numcases}{}
 -\frac{1}{2}A^{0}_{001,001}+\frac{1}{2}A^{0}_{010,010}+\frac{1}{4}A^{2}_{011,010}=2b \label{g}\\
 -\frac{1}{2}A^{0}_{001,001}+\frac{1}{4} A^{1}_{011,001}+\frac{1}{2}A^{0}_{010,010}=0  \label{h}
\end{numcases}

 On the other hand, for $(g,h)=(\theta_2,\theta_2)$ and for $(g,h)=(\theta_1,\theta_1)$ we consider the terms in ${\partial_x f}$ in \eqref{eq1} then we find
\begin{numcases}{}
\frac{1}{4}A^{0}_{000,100}+\frac{1}{2}A^{2}_{001,100}-\frac{3}{2}A^{0}_{001,001}=-4b \label{i}. \\
\frac{1}{4}A^{0}_{000,100}+\frac{1}{2}A^{2}_{010,100}-\frac{3}{2}A^{0}_{010,010}=0 \label{j}.
\end{numcases}
 For $(g,h)=(\theta_2,\theta_1 \theta_2)$ the terms in  ${\theta_1\partial_x f}$ give
 \begin{align}\label{k}
  -\frac{3}{4}A^{2}_{001,100}-\frac{1}{4}A^{12}_{011,100}-\frac{1}{2}A^{0}_{000,100}+\frac{3}{4}A^{0}_{001,001}-\frac{3}{4}A^{1}_{011,001}=2b
 \end{align}
 and for  $(g,h)=(\theta_1,\theta_1 \theta_2)$ the terms in ${\theta_2\partial_x f}$ imply
\begin{align}\label{l}
\frac{3}{4}A^{2}_{010,100}+\frac{1}{4}A^{12}_{011,100}+\frac{1}{2}A^{0}_{000,100}-\frac{3}{4}A^{0}_{010,010}-\frac{3}{4}A^{2}_{011,010}=0
\end{align}
Now, combining equations coming from substituting \eqref{a} into \eqref{b}, adding \eqref{e} and \eqref{f},  \eqref{g} and \eqref{h}, substituting \eqref{i} into \eqref{j} and adding \eqref{k} and \eqref{l}, we immediately find  $b=0$ \\

To complete the proof we  proceed similarly as before, therefore we get
\begin{numcases}{}
\frac{1}{4}A^{1}_{101,000}+\frac{1}{4}A^{2}_{110,000}=a \nonumber\\
\frac{1}{4}A^{1}_{101,000}+\frac{1}{2}A^{2}_{110,000}=a \nonumber\\
\frac{1}{2}A^{1}_{101,000}+\frac{1}{4}A^{2}_{110,000}=a. \nonumber
 \end{numcases}
Thus, it is easy to see that $a=0$. So we obtain the claim.  \qed
\end{proofname}
\subsection{Case 2: $2d=m \in \mathbb{N}$}
\begin{theorem} \label{th2}
If $2d=m \in \mathbb{N}$ then the image $\mathrm{H}^1\vee\mathrm{H}^1$ of $\mathrm{H}^{1}_{\mathrm{diff}}\left(\mathfrak{osp}(2|2), \mathfrak{D}_{d-\frac{k}{2},d-\frac{k}{2}}\right)$ under the cup-product is a 6-dimensional subspace of $\mathrm{H}^{2}_{\mathrm{diff}}\left(\mathfrak{osp}(2|2), \mathfrak{D}_{-\frac{k}{2},\frac{k}{2}}\right)$ spanned by
$$\Phi_{1}=\gamma_{k}\vee \widetilde{\Gamma}_{k},\; \Phi_2=\gamma_{k}\vee \overline{\Gamma}_{k}\;,\Phi_3=  \widetilde{\gamma}_{k}\vee\widetilde{\Gamma}_{k},\;\Phi_4=\widetilde{\gamma}_{k}\vee\overline{\Gamma}_{k},\; \Phi_5= \widetilde{\Gamma_{k}}\vee \gamma_{-k}
\;\mbox{and}\;\Phi_6=\overline{\Gamma}_{k}\vee \widetilde{\gamma}_{-k}.$$
 \end{theorem}
\begin{proofname}. In this case, the space $\mathrm{H}^1\vee\mathrm{H}^1$ is generated by the following twelve cup-products:\\ $\Phi_{1}=\gamma_{k}\vee \widetilde{\Gamma}_{k}$, $\Phi_2=\gamma_{k}\vee\overline{\Gamma}_{k}$, $\Phi_3=  \widetilde{\gamma}_{k}\vee\widetilde{\Gamma}_{k}$, $\Phi_4=\widetilde{\gamma}_{k}\vee\overline{\Gamma}_{k}$, $\Phi_5= \widetilde{\Gamma}_{k}\vee\gamma_{-k},$ $\Phi_6= \overline{\Gamma}_{k}\vee \widetilde{\gamma}_{-k}$, $\Phi_7=\gamma_{k}\vee \Gamma_{k}$, $\Phi_8= \Gamma_{k}\vee \widetilde{\gamma}_{-k}$, $\Phi_9=  \widetilde{\Gamma}_{k}\vee \widetilde{\gamma}_{-k}$, $\Phi_{10}=  \Gamma_{k}\vee \gamma_{-k}$, $\Phi_{11}=  \overline{\Gamma}_{k}\vee \gamma_{-k}$ and
  $\Phi_{12}=\widetilde{\gamma}_{k}\vee\Gamma_{k}$.

  By a straightforward computation, we check that
 $$\Phi_7=0,\quad     \Phi_8=-\Phi_1,\quad   \Phi_{10}=\Phi_{11}=-\Phi_2,\quad    \Phi_9=-\Phi_3,\quad\text{and}\quad
  \Phi_{12}=\Phi_4+\Phi_5$$
    where
   {\small\begin{eqnarray*}
          \Phi_{1}(g,h) &=&(-1)^{k} \Big[(-1)^{h}g'\partial_{2}h\theta_{2}-(k+1)g'\overline{\eta}_{1}\partial_{2}h\theta_{1}\theta_{2}\Big]
           \mathbf{\partial_{x}^{k+1}} +(-1)^{k} \textbf{[}kg'\overline{\eta}_{1}\partial _2h\theta_1-(-1)^{h}g'\partial_2H\textbf{]}\mathbf{\partial^{k}_x\partial_2}\\
             && -(-1)^{k}(k+1)g' \overline{\eta}_{1}\partial _2h\theta_2\mathbf{\partial^{k}_x\partial_1}-(-1)^{k}kg'\overline{\eta}_{1}\partial _2H\mathbf{\partial_{x}^{k-1}\partial_1\partial_2}-(-1)^{gh}(g\leftrightarrow h)
          \end{eqnarray*}}
         {\small\begin{eqnarray*}
           \Phi_2(g,h)&=& (-1)^k\Big[(k-1)g'h''\theta_2-(-1)^{h}g'\overline{\eta}_{2}h'\Big]\mathbf{\partial^{k-1}_x\partial_1}  -(-1)^k\Big[(k-1)g'h''\theta_1-(-1)^{h}g'\overline{\eta}_{1}h'\Big]\mathbf{\partial^{k-1}_x\partial_2} \\
            && +(-1)^k\Big[(k-1)g'h''\theta_1\theta_2+(-1)^{h}g'\overline{\eta}_{2}h'\theta_1-(-1)^{h}g'\overline{\eta}_{1}h'\theta_2\Big] \mathbf{\partial_{x}^{k}} \\
            &&+ (-1)^k(k-1)g'h''\mathbf{ \partial^{k-2}_x\partial_1 \partial_2} -(-1)^{gh}(g\leftrightarrow h)
         \end{eqnarray*}}
       {\small\begin{eqnarray*}
          \Phi_{3}(g,h) &=&(-1)^k\Big[(-1)^{g+h}\partial_2g\overline{\eta}_{1}\partial_{2}h+(-1)^{g}k\partial_2g\partial_{x}\partial_2h
            \theta_1\Big]\mathbf{\partial_x^{k}\partial_2} -(-1)^{g+k}(k+1)  \partial_2g\partial_x\partial_2h\theta_2\mathbf{\partial_x^{k}\partial_1}\\&&-(-1)^k\Big[(-1)^{g+h}\partial_2g\partial_{1}\partial_{2}h\theta_2
           +(-1)^{g}(k+2)\partial_2g\partial_{x}\partial_2h\theta_1\theta_2\Big]\mathbf{\partial_x^{k+1}}   -(-1)^{k+g} \partial_2g\partial_2h\times\\
           &&\Big[ \mathbf{\partial_1\partial_2}+\theta_2\mathbf{\partial_{x}\partial_1}
           -\theta_1\mathbf{\partial_{x}^{k+1}\partial_2}+\theta_1\theta_2 \mathbf{\partial^{2}_x } \Big]  -(-1)^{k+g}k\partial_2g\partial_x\partial_2h\mathbf{\partial_x^{k-1}\partial_1 \partial_2}
             -(-1)^{gh}(g\leftrightarrow h)
\end{eqnarray*}}
{\small\begin{eqnarray*}
 \Phi_{4}(g,h)&=&  (-1)^{k}\Big[(-1)^{g}\partial_2g\overline{\eta}_{1}h'\Big]\Big[\theta_2\mathbf{\partial_x^{k}\partial_1} +\theta_1\theta_2\mathbf{\partial_x^{k+1}}\Big] -(-1)^{k}\Big[k(k-1)\overline{\eta}_{1}\partial_2gh''\theta_1-(-1)^{h}k\overline{\eta}_{1}\partial_2g\overline{\eta}_{1}h'\\
   &&
  -(-1)^{g+h}k\partial_2gh''\Big]\mathbf{\partial_x^{k-1}\partial_2}  +(-1)^{k}\Big[k^2\overline{\eta}_{1}\partial_2gh''\theta_2
   -(-1)^{h}k\overline{\eta}_{1}\partial_2g\overline{\eta}_{2}h'\Big]\mathbf{\partial_x^{k-1}\partial_1} \\
   &&+(-1)^{k}\Big[k^2\overline{\eta}_{1}\partial_2gh''\theta_1\theta_2
   +(-1)^{h}k\overline{\eta}_{1}\partial_2g\overline{\eta}_{2}h'\theta_1
   -(-1)^{h}(k+1)\overline{\eta}_{1}\partial_2g\overline{\eta}_{1}h'\theta_2  \\
   && -(-1)^{g+h}k\partial_2gh''\theta_2+(-1)^{g}\partial_2g\overline{\eta}_{2}h'\Big]\mathbf{\partial_x^{k}}+(-1)^{k}k(k-1)\overline{\eta}_{1}\partial_2gh'' \mathbf{\partial_x^{k-2}\partial_1\partial_2}\\
  &&+
   (-1)^{k}\Big[(-1)^{g}\partial_2(g)\overline{\eta}_{1}h'+(-1)^{g}\theta_2\overline{\eta}_{2}\overline{\eta}_{1}g\overline{\eta}_{2}h'\textbf{]}
 \Big[\mathbf{\partial_x^{k-1}\partial_1\partial_2}-\theta_1\mathbf{\partial_x^{k}\partial_2} \Big]-(-1)^{gh}(g\leftrightarrow h).
\end{eqnarray*}}
{\small\begin{eqnarray*}
  \Phi_5(g,h)&=& (-1)^{k}\Big[k(k-1)\overline{\eta}_{1}\partial_{2}gh''\theta_1 -(-1)^{h}k\overline{\eta}_{1}\partial_{2}g\overline{\eta}_{1}g'-(-1)^{g+h}k\partial_{2}gh''
   \Big]\mathbf{\partial^{k-1}_x\partial_2}   -(-1)^{k}\Big[k^2\overline{\eta}_{1}\partial_{2}gh''\theta_2\\
   && -(-1)^{h}k \overline{\eta}_{1}(
   \partial_{2}(g)\overline{\eta}_{2}h'\Big]\mathbf{ \partial^{k-1}_x\partial_1} -(-1)^{k}k(k-1)\overline{\eta}_{1}\partial_2gh''\mathbf{\partial^{k-2}_x\partial_1 \partial_2}-(-1)^{k}k\overline{\eta}_{1}\partial_2gh'\mathbf{\partial^{k-1}_x\partial_1 \partial_2}\\
   && -(-1)^{k} \Big[k\overline{\eta}_{1}\partial_{2}gh'\theta_2-(-1)^{g+h}\overline{\eta}_{1}   \theta_2\partial_{2}gh'\Big]   \mathbf{\partial^{k}_x\partial_1}+ (-1)^{k}\Big[k\overline{\eta}_{1}\partial_{2}gh'\theta_1-(-1)^{h+g}\partial_{2}gh' \Big] \mathbf{ \partial^{k}_x\partial_2}   \\
   &&-(-1)^{k}\Big[k\overline{\eta}_{1}\partial_{2}gh'\theta_1\theta_2-(-1)^{g+h}\partial_{2}(g)h'\theta_2
   +(-1)^{g+h}\overline{\eta}_{1}(\theta_2\partial_{2}g)h'\theta_1\Big] \mathbf{\partial^{k+1}_x}  \\
   &&+(-1)^{k}\Big[(-1)^{h}(k+1)\overline{\eta}_{1}\partial_{2}g\overline{\eta}_{1}g'\theta_2
   -k^2\overline{\eta}_{1}\partial_{2}gh''\theta_1\theta_2+(-1)^{g+h}k\partial_{2}gh''\theta_2 - (-1)^{h}k \overline{\eta}_{1}\partial_{2}g\overline{\eta}_{2}h'\theta_1\\
   &&-(-1)^{g}\partial_2g\overline{\eta}_{2}h'
   \Big]\mathbf{\partial^{k}_x} - (-1)^{gh}(g\leftrightarrow h)
\end{eqnarray*}}
 {\small \begin{eqnarray*}
 \Phi_{6}(g,h) &=& -(-1)^k\Big[(-1)^{g}\overline{\eta}_{2}g'\theta_2\overline{\eta}_{2}\overline{\eta}_{1}h +(-1)^{g}\overline{\eta}_{1}g'\partial_2h\Big]
 \Big[\mathbf{\partial_x^{k-1}\partial_1 \partial_2}+\theta_2\mathbf{\partial_x^{k}\partial_1}-\theta_1\mathbf{\partial_x^{k}\partial_2}
 +\theta_1\theta_2 \mathbf{\partial_x^{k+1}}\Big]\\&&-(-1)^{k}\Big[k(k-1)g''\overline{\eta}_{1}\partial_{2}h+(-1)^{g}(k-1)\overline{\eta}_{1}g'\partial_{x}\partial_{2}h\Big] \mathbf{ \partial_x^{k-2}\partial_1 \partial_2} -\Big[(k-1)(k+1)g''\overline{\eta}_{1}\partial_2h\theta_2\\
 &&- (-1)^{g+h}k\overline{\eta}_{1}g'\theta_2\partial_x\partial_2h  -(-1)^{g+h}(k+1)\overline{\eta}_{2}g'\overline{\eta}_{1}\partial_2h\Big](-1)^{k}\mathbf{\partial_x^{k-1}\partial_1}  +\Big[k(k-1)g''\overline{\eta}_{1}\partial_{2}h\theta_1\\
 &&
 +(-1)^{g}(k-1)\overline{\eta}_{1}g'\partial_x\partial_2h\theta_1-(-1)^{h}(k-1)g''\partial_2h-(-1)^{g+h}
 \overline{\eta}_{2}g'\theta_2\partial_2\partial_x \\
 &&-(-1)^{g+h}(k-1)\overline{\eta}_{1}g' \overline{\eta}_{1}\partial_2h
\Big](-1)^{k}\mathbf{\partial_x^{k-1}\partial_2}-(-1)^k\Big[(k+1)(k-1)g''\overline{\eta}_{1}\partial_{2}(h)\theta_1\theta_2\\&&
 -(-1)^{g+h}k\overline{\eta}_{1}g'\overline{\eta}_{1}\partial_{2}h\theta_2
 +(-1)^{g}k\overline{\eta}_{1}g'\partial_2\partial_x h\theta_1\theta_2
 -(-1)^{h}(k-1)g''\partial_2h\theta_2\\
 &&
 +(-1)^{g+h}(k+1)\overline{\eta}_{2}g'\partial_1\partial_2h\theta_1
 -(-1)^{g}\overline{\eta}_{2}g'\partial_2h\Big]\mathbf{\partial_x^{k}} +(-1)^{h+k}(k-1)g''\partial_x\partial_2h \mathbf{ \partial_x^{k-2}\partial_2}\\&& + (-1)^{k}\Big[-(-1)^{h}(k-1)g''\partial_x\partial_2h\theta_2-
  (-1)^{g}\overline{\eta}_{2}g'\partial_x\partial_2h\Big]\mathbf{\partial_x^{k-1}} -(-1)^{gh}(g\leftrightarrow h)
  \end{eqnarray*}}

The proof is almost identical to the previous theorem. Indeed, we have to prove that  $\Phi_1$, $\Phi_2$, $\Phi_3$, $\Phi_4$, $\Phi_5$, and $\Phi_{6}$ are nontrivial 2-cocycles which are linearly independent. That is, the equation
\begin{equation}\label{eq2}
a_1\Phi_1+a_2\Phi_2+a_3\Phi_3+a_4\Phi_4+a_5\Phi_5+a_6\Phi_{6}=\delta B
\end{equation}
 has a solution if and only if $a_1=\dots=a_6=0$. Of course we can express $B$ as in \eqref{cobo}.

Here, we just give the parameters of identifications which allow us to obtain the result.
For \eqref{eq2}, considering the terms in $\partial_{x}^{k}\partial_{2}$ for
$(g,h)=(1,\theta_2)$ and  $(g,h)=(x,\theta_2)$ and the terms in $\partial_{x}^{k}\partial_{1}$ for
$(g,h)=(1,\theta_1)$ and $(g,h)=(x,\theta_1)$, we get  $A^{2}_{100,k01}=-4(-1)^{k}(a_1-a_5)$.
   But the terms in $ \theta_1\partial_{x}^{k}\partial_{2}$, for $(g,h)=(1,\theta_1\theta_2)$ and then for $(g,h)=(x,\theta_1\theta_2)$ give $A^{2}_{100,k01}=(-1)^{k}2(k-1)(a_1-a_5)$, therefore $a_1=a_5$. \\
 Similarly, considering  the terms in $\theta_1\theta_2\partial_{x}^{k}$, for
$(g,h)=(1,x), \, (1,x^2),\,(x,x^2)$ we obtain  $a_2=0$.

 Now, since $a_1=a_5$ and $a_2=0$, the terms in $\theta_1\theta_2\partial^{k}_x$, for
   $(g,h)=(1,\theta_1\theta_2), \,(x,\theta_1\theta_2),\,(x^2,\theta_1\theta_2)$ give $a_4-a_6+a_5=0$. Thus, the terms in $\partial_{x}^{k-2}\partial_1\partial_2$, for $(g,h)=(x\theta_1,x\theta_2),\,(\theta_1,\theta_2),\,(x\theta_1,\theta_2),\,(x\theta_2,\theta_1)$
give $A^{0}_{101,(k-1)01}-A^{0}_{110,(k-1)10}=(-1)^{k}4a_6$.
On the other hand, the terms in $\partial_{x}^{k-1}\partial_{2}$ for $(g,h)=(\theta_1, \theta_1\theta_2),\,(x\theta_1, \theta_1\theta_2)$ and the terms in $\partial_{x}^{k-1}\partial_{1}$ for  $(g,h)=(\theta_2, \theta_1\theta_2),\,(x\theta_2, \theta_1\theta_2)$, give
 $A^{0}_{101,(k-1)01}-A^{0}_{110,(k-1)10}=-(-1)^{k}\frac{4}{3}a_6$, then $a_6=0.$ So, considering the terms in  $\partial_{x}^{k}\partial_1\partial_2$, for $(g,h)=(\theta_1,\theta_2)$  combined with those in $\theta_1\partial_{x}^{k-1}\partial_1\partial_2$, for $(g,h)=(x^2,\theta_2),\,(x^2,\theta_1)$ and  by those in $\theta_1\partial_{x}^{k-1}\partial_1\partial_2$, for $(g,h)=(x\theta_1,\theta_2),\,(x\theta_2,\theta_1)$, we get $$A^{1}_{101,(k-1)11}+A^{2}_{110,(k-1)11}= -(-1)^k2a_4+k (A^{0}_{010,k10}-A^{0}_{001,k01}) $$  and $$A^{1}_{101,(k-1)11}+A^{2}_{110,(k-1)11}=-(-1)^k4a_4
+k (A^{0}_{010,k10}-A^{0}_{001,k01}).$$
 Then $a_4=0$ and consequently $a_1=a_5=0$.

Finally, we consider the terms in $\partial_{x}^{k}\partial_2$, for $(g,\,h)=(\theta_2,\theta_1\theta_2)$, the terms in $\partial_{x}^{k}\partial_1$, for $(g,\,h)=(\theta_1,\theta_1\theta_2)$, the terms in $\partial_{x}^{k-1}\partial_1\partial_2$, for $(g,\,h)=(x\theta_2,\theta_1\theta_2),\,(x\theta_1,\theta_1\theta_2),\,(x\theta_2,\theta_2),\,(x\theta_1,\theta_1)$ and the terms in $\partial_{x}^{k}\partial_1\partial_2$, for
 $(g,\,h)=(\theta_2,\theta_1\theta_2),\,(\theta_1,\theta_1\theta_2),\,(\theta_2,\theta_2),\,(\theta_1,\theta_1)$, then we have $$A^{0}_{001,k10}+A^{0}_{010,k01}=-(-1)^{k}\frac{2}{3}a_3$$ and  $$A^{0}_{010,k01}+A^{0}_{001,k10}=-(-1)^{k}a_3.$$ Then $a_3=0.$ \qed
 \end{proofname}

 We formulate a conjecture on the structure of the second cohomology space.
\begin{conjecture} One has $\mathrm{H}^1_{\mathrm{diff}}\left(\mathfrak{osp}(2|2), \mathfrak{D}_{\lambda,\mu}\right)\vee\mathrm{H}^1_{\mathrm{diff}}\left(\mathfrak{osp}(2|2), \mathfrak{D}_{\lambda,\mu}\right)=\mathrm{H}^{2}_{\mathrm{diff}}\left(\mathfrak{osp}(2|2), \mathfrak{D}_{\lambda,\mu}\right)$.
\end{conjecture}
 This conjecture is an important open problem concerning the computation of second cohomology spaces which are generally difficult to derive. It turns out that a positive confirmation of this type of conjecture is a crucial result as obtained by Arnal, Ben Ammar and Dali  \cite{ddar} where they proved that $\mathrm{H}^2(\mathfrak{osp}(1|2),\,\mathfrak{D}_{\lambda,\mu})=\mathrm{H}^1(\mathfrak{osp}(1|2),\,\mathfrak{D}_{\lambda,\mu})\vee \mathrm{H}^1(\mathfrak{osp}(1|2),\,\mathfrak{D}_{\lambda,\mu})$.

\section{Integrability Conditions}
Now, we consider an infinitesimal deformation
\begin{equation}\label{inf1}\widetilde{\mathfrak{L}}=\mathfrak{L}+ \mathfrak{L}^{1} \end{equation}
of the natural action $\mathfrak{L}$ of $\mathfrak{osp}(2|2)$ on the space
 $\mathfrak{S}_d^{2}$ and we study the necessary and sufficient conditions to extend it to a formal one:
\begin{equation}
\label{BigDef3} \widetilde{\mathfrak{L}}= \mathfrak{L}+\mathfrak{L}^{1}+\sum_iP_i^{2}\mathfrak{L}_i^{2}+ \sum_i P_i^{3}\mathfrak{L}_i^{3}+\cdots
\end{equation}
 where
 $$\mathfrak{L}^{1}= \left\{
\begin{array}{ll}
\sum_{k\geq 0} (a_{k}\omega_{k}+b_{k}\widetilde{\omega}_{k})&\text{if}\quad 2d\notin\mathbb{N}\\
\sum _{k\leq
m}(a_{k}\gamma_{k}+b_{k}\widetilde{\gamma}_{k})+\sum
_{k=1}^{m}(c_{k}\Gamma_{k}+d_{k}\widetilde{\Gamma}_{k}+e_{k}\overline{\Gamma}_k) &
\text{if}\quad
2d=m \in \mathbb{N}
 \end{array}
 \right.$$
 and the the higher order terms
$\mathfrak{L}_i^{2}$, $\mathfrak{L}_i^{3},\ldots$ are linear maps from
$\mathfrak{osp}(2|2)$ to $\mathrm{End}(\mathfrak{S}_d^{2})$ such that  the map
\begin{equation} \label{map1} \widetilde{\mathfrak{L}}:\mathfrak{osp}(2|2)\to
\mathbb{C}[[a_{k},b_{k},c_{k},d_{k},e_{k}]]\otimes{\rm End(\mathfrak{S}_d^{2})},
\end{equation}
satisfies the homomorphism condition: \begin{equation}\label{hom}\widetilde{\mathfrak{L}}_{[f,g]}= [\widetilde{\mathfrak{L}}_f,\widetilde{\mathfrak{L}}_g],\end{equation} and $P_i^{j}$ are monomial in the independent parameters $a_{k},b_{k},c_{k},d_{k},e_{k}$ (or $a_{k},b_{k}$ if $2d\notin\mathbb{N}$) with degree $j$ and with the same parity of $\mathfrak{L}_i^{j}$.\\
The following theorems are our main results.  We have to distinguish two cases.

\subsection{Case 1: $2d \notin \mathbb{N}$}
In this case, we have  $\mathfrak{L}^{1}=
\sum_{k\geq 0} (a_{k}\omega_{k}+b_{k}\widetilde{\omega}_{k})$.
\begin{theorem}
The following conditions are necessary and sufficient  for integrability of the infinitesimal deformation  \eqref{inf}
\begin{equation}\label{cond1}
b_k=0,\quad\text{for all}\quad k\geq0.
 \end{equation}
 Moreover, any formal deformation is equivalent to its infinitesimal part which is of the form:
  \begin{equation}\label{thmcase1}\mathfrak{L}+\sum_{k\geq 0} a_{k}\omega_{k}.\end{equation}
\end{theorem}
That is, formal deformations are classified by the subspace of $\mathrm{H^1}(\mathfrak{osp}(2|2),\frak{D}^2_{\lambda,\mu})$ spanned by the cohomological classes of the 1-cocycles  $\omega_{k}$.\\

\begin{proofname}.
The condition \eqref{hom} gives, for the second-order terms,
the following equation
\begin{equation}\label{ttttt}
\delta{\mathfrak{L}^{2}}=\frac{1}{2}\sum_{k\geq 0}(a_kb_k\Omega_1+b^{2}_k\Omega_2).\end{equation}
\end{proofname}
Thus, the right hand side of \eqref{ttttt} must be
a coboundary. But, by Theorem \ref{th1},  $\Omega_1$ and $\Omega_2$ are linearly independent nontrivial 2-cocycles, therefore $a_kb_k=b^{2}_k=0$ for all $k\geq0$. Thus, the conditions \eqref{cond1} are necessary.

Now, we show that these conditions are sufficient. The solutions
$\frak{L}^{k}$ of the Maurer-Cartan equations \eqref{maurrer cartank} are defined up to a
1-cocycle and it has been shown in works \cite{aalo} and \cite{ff} that different
choices of solutions correspond to
equivalent deformations. Thus,  we can always reduce
 $\mathfrak{L}^{k}$, for $k=2$ to zero by equivalence. Then, by recurrence, the highest-order terms $\mathfrak{L}^{k}$ with $k\geq 3$, also satisfy the equation $\delta(\mathfrak{L}^{k})$ and can also be reduced to the identically zero map. One obviously
obtains a deformation (which is of order 1 in $a_k$).\qed
 \subsection{Case 2: $2d=m \in \mathbb{N}$}
  In this case we have $\mathfrak{L}^{1}= \sum _{k\leq
m}(a_{k}\gamma_{k}+b_{k}\widetilde{\gamma}_{k})+\sum
_{k=1}^{m}(c_{k}\Gamma_{k}+d_{k}\widetilde{\Gamma}_{k}+e_{k}\overline{\Gamma}_k)$.
 \begin{theorem}
The following conditions are necessary and sufficient  for integrability of the infinitesimal deformation  \eqref{inf}: For all $k\geq0$,
\begin{numcases}{}
a_k d_k -c_k b_{-k}=0\nonumber \\
a_k e_k-c_k a_{-k}-e_k a_{-k}=0\nonumber \\
b_k d_k-d_k b_{-k} =0 \nonumber\\
b_k c_k=- b_k e_k=-d_k a_{-k}\nonumber\\
 e_k b_{-k}=0. \nonumber
\end{numcases}
Moreover, any formal deformation is equivalent to its infinitesimal part.
\end{theorem}
Thus, similarly to the first case, these above conditions give a classification of the formal deformations.\\

 \begin{proofname}.
 In this case, the equation \eqref{hom} can be expressed as follows
 \begin{eqnarray}\label{vvvvv}
  \delta{\mathfrak{L}^{2}}&=&\frac{1}{2}\sum_{k\geq 0}\big[(a_k d_k -c_k b_{-k})\Phi_1+(a_k e_k-c_k a_{-k}-e_k a_{-k})\Phi_2+ (b_k d_k-d_k b_{-k})\Phi_3\nonumber\\
  &+&(b_k c_k+ b_k e_k)\Phi_4
  +( b_k c_k+d_k a_{-k})\Phi_5+ e_k b_{-k}\Phi_6\big]
 \end{eqnarray}
 The second order integrability conditions are determined
by the fact that the map 2-cocycles  $\Phi_1$, $\Phi_2$, $\Phi_3$, $\Phi_4,$ $\Phi_5$, $\Phi_6$ are non-trivial, which is proved in Theorem \ref{th2}. As above, these conditions are sufficient and the terms $\mathfrak{L}^{k}$ with $k\geq2$ can be chosen identically
zero.
 \end{proofname}
\begin{example}
Let us consider $d=\frac{m}{2}\in \mathbb{N}$, let $b_k=d_k=0$; $k\in\mathbb{Z}$ and $c_{k}=e_{k}$. So, we obtain the following deformation of  $\mathfrak{S}_d^{2}$ with two family of independent parameters
\begin{equation*}
  \widetilde{\mathfrak{L}}=\mathfrak{L}+ \sum _{k\leq
m}2 a_{-k}\gamma_{k}+\sum
_{k=1}^{m}c_{k}(\Gamma_{k}+\overline{\Gamma}_k).
\end{equation*}
 Of course it is easy to give many other examples of true deformations.
\end{example}

\bigskip\noindent



\begin{thebibliography}{20}
\bibitem{akl}
 M.Abdaoui, H. Khalfoun and I.Laraeidh,
\newblock Deformation of modules of weighted densities on the superespace $\mathbb{R}^{1|N}$,
\newblock {\em Acta Mathematica Hungarica.}
 \textbf{145}(2015), no. 1, 104-123.

\bibitem{aalo}
B. Agrebaoui, F. Ammar, P. Lecomte, V. Ovsienko,
\newblock Multi-parameter deformations of the module of symbols of
differential operators,
\newblock {\em Internat. Mathem. Research Notices.}
 \textbf{16} (2002), 847--869.

\bibitem{abbo}
B. Agrebaoui, N. Ben Fraj, M. Ben Ammar, and V. Ovsienko,
\newblock Deformation of modules of differential forms,
\newblock {\em NonLinear Mathematical Physics.}
 \textbf{ 10} (2003), no. 2, 148--156.

\bibitem{ddar}
D. Arnal, M. Ben Ammar, and Bechir Dali,
 \newblock The spaces $\mathrm{H}^n(\mathfrak{osp}(1|2),M)$ for some weight modules $M$,
\newblock {\em J. Math. Phys.}
\textbf{ 51}(2010).

 \bibitem{bb1}
 I. Basdouri, M. Ben Ammar,
 \newblock Cohomology of $\frak {osp}(1|2)$
Acting on Linear Differential Operators on the Supercercle $S^{1|1}$,
\newblock {\em Letters in Mathematical Physics.}
\textbf{81}(2007), 239--251.

\bibitem{bb2}
I. Basdouri, M. Ben Ammar,
\newblock Deformation of $\mathfrak{sl}(2)$ and $\mathfrak{osp}(1|2)$-Modules of Symbols,
\newblock {\em Acta Math. Hungar.}
\textbf{ 137}(2012), no 3, 214-223.

\bibitem{bbbbk}
I. Basdouri, M. Ben Ammar, N. Ben Fraj, M. Boujelbene and K. Kammoun,
 \newblock Cohomology of the Lie Superalgebra of Contact Vector
Fields on $\mathbb{R}^{1|1} $ and Deformations of the Superspace of Symbols,
\newblock {\em Journal of Nonlinear Math Physics.}
\textbf{ 16}(2009), No. 4, 1--37.

\bibitem{bbdo}
I. Basdouri, M. Ben Ammar, B. Dali and S. Omri,
\newblock  Deformation of $\mathfrak{vect}(1)$-Modules of Symbols,
\newblock {\em Journal of Geometry and Physics.}
\textbf{60}(2010), 531--542.

\bibitem{bb}
M. Ben Ammar, M. Boujelbene,
\newblock $\mathfrak{sl}(2)$-Trivial Deformation of $\mathrm{Vect_{Pol}}(\mathbb{R})$-Modules of Symbols, SIGMA
\textbf{ 4}(2008), 065-- 19 pages.

\bibitem{bb3}
N. Ben Fraj, M. Boujelben,
\newblock Cohomology of $\mathfrak{osp}$(2|2) Acting on Spaces of Linear Differential Operators on the Superspace $\mathbb{R}^{1|2}$,
 ISSN 0001-4346, Mathematical Notes, 2012, Vol. 92, No. 3, pp. 302-311.  Pleiades Publishing, Ltd., 2012.Published in Russian in Matematicheskie Zametki, 2012, Vol. 92, No. 3, pp. 331-342.

 \bibitem{ff}
 A. Fialowski, D. B. Fuchs,
\newblock Construction of miniversal deformations of Lie algebras,
\newblock{\em J. Func. Anal.} {\bf 161:1} (1999), 76--110.

\bibitem{gro}P. Grozman, D. Leites, and I. Shchepochkina,
\newblock Lie superalgebras of string theories,
\newblock{\em ActaMath. Vietnam} 26
(1), 27–63 (2001).

\bibitem{lec}P. B. A. Lecomte, {\it On the cohomology of
$\frak{sl}(n + 1;\mathbb{R})$ acting on differential operators and
$\frak{sl}(n + 1;\mathbb{R})$-equivariant symbols,} Indag. Math. NS.
11 (1), (2000), 95 114.

\bibitem{nr}
A. Nijenuis, R. W. Richardson Jr., {\it Deformations of
homomorphisms of Lie groups and Lie algebras}, Bull. Amer. Math.
Soc. {\bf 73} (1967), 175--179.

\end{thebibliography}
\end{document}